\documentclass[12pt,a4paper]{article}
\usepackage[T1]{fontenc}
\usepackage{mathptmx}
\usepackage{courier}
\usepackage[dvips]{graphicx}
\usepackage{psfrag}
\usepackage{amsmath,amssymb}        
\usepackage{latexsym}
\usepackage{url}

\def\ds{\displaystyle}

\usepackage[T1]{fontenc}
\usepackage{mathptmx}
\usepackage{courier}
\usepackage[dvips]{graphicx}
\usepackage{psfrag}
\usepackage{amsmath,amssymb}        
\usepackage{latexsym}
\usepackage{url}

\def\ds{\displaystyle}

\title{Numerical Transformation Methods\\ for a Moving-Wall Boundary Layer Flow of a\\
Rarefied Gas Free Stream\\ over a Moving Flat Plate}
\author{Riccardo Fazio \\
Department of Mathematics, Computer Science\\ Physical Sciences and Earth Sciences,\\
University of Messina \\
Viale F. Stagno D'Alcontres, 31 \\
98166 Messina, Italy \\
E-mail: rfazio@unime.it \\
Home-page: http://mat521.unime.it/fazio}
\date{\today}
\begin{document}
\maketitle
\begin{abstract}
The first contribution of this paper is the extension of the non-iterative transformation method, proposed by T\"opfer more than a century ago and defined for the numerical solution of the Blasius problem, to a Blasius problem with extended boundary conditions.
This method, which makes use of the invariance of two physical parameters with respect to a scaling group of point transformation, allows us to solve numerically the Blasius problem with extended boundary conditions by solving a related initial value problem and then rescaling the obtained numerical solution. 
Therefore, our method is an initial value method.
However, in this way, we cannot fix in advance the physical parameters, and if we need just to compute the numerical solution for given values of the two parameters we have to define an iterative extension of the transformation method, which is the second contribution of this work.
\end{abstract}
\bigskip

\bigskip


\section{Introduction}
It was Prandtl \cite{Prandtl:1904:UFK} in 1904 who fixed the terms of validity of boundary layer theory.
Within this theory, the problem of determining the steady two-dimensional motion of fluid flow past a flat plate placed edge-ways to the mainstream was
formulated in general terms and investigated in detail by Blasius \cite{Blasius:1908:GFK}. 
The engineering interest was to calculate the shear at the plate (skin friction), which leads to the determination of the viscous drag on the plate, see for instance Schlichting \cite{Schlichting:2000:BLT}).
Blasius problem is a boundary value problem (BVP) defined on the semi-infinite interval $[0, \infty)$.
It is possible to prove, see Weyl \cite{Weyl:1942:DES}, that the unique solution of the Blasius problem has a positive second-order derivative, which is a monotone decreasing function on $ [0, \infty) $ and approaches to zero as $ \eta $ goes to infinity.
The governing differential equation and the two boundary conditions at the origin in the Blasius problem are invariant with respect to a scaling group of transformations and this has several consequences.
From a numerical viewpoint a non-iterative transformation method (ITM)
reducing the solution of the Blasius problem to the solution of a related initial value problem (IVP) was defined by T\"opfer \cite{Topfer:1912:BAB}. 
Consequently, applying the scaling invariance properties, a simple existence and uniqueness Theorem
was given by J. Serrin, see Meyer \cite[pp. 104-105]{Meyer:1971:IMF}. 
Let us note here that the mentioned invariance property is essential to the error analysis of the truncated boundary solution due to Rubel \cite{Rubel:1955:EET}, see Fazio \cite{Fazio:2002:SFB} for the full details.
Recently, the Blasius problem was used, by Boyd \cite{Boyd:2008:BFC}, as an example, were some good analysis allowed researchers of the past to solve problems, governed by partial differential equations, that might be otherwise impossible to face before the computer invention.

Non-ITMs have been applied to several problems of practical interest within the applied sciences. 
First of all, a non-ITM was applied to the Blasius equation with slip boundary condition, arising within the study of gas and liquid flows at the micro-scale regime \cite{Gad-el-Hak:1999:FMM,Martin:2001:BBL}, see \cite{Fazio:2009:NTM}.
A non-ITM was applied also to the Blasius equation with moving wall considered by Ishak et al. \cite{Ishak:2007:BLM} or surface gasification studied by Emmons \cite{Emmons:1956:FCL} and recently by Lu and Law \cite{Lu:2014:ISB} or slip boundary conditions investigated by Gad-el-Hak \cite{Gad-el-Hak:1999:FMM} or Martin and Boyd \cite{Martin:2001:BBL}, see Fazio \cite{Fazio:2016:NIT} for details.
In particular, within these applications, we found a way to solve non-iteratively the Sakiadis problem \cite{Sakiadis:1961:BLBa,Sakiadis:1961:BLBb}.
The application of a non-ITM to an extended Blasius problem has been the subject of a recent paper \cite{Fazio:2020:NIT}.
As far as the non-ITM is concerned, a recent review dealing with all the cited problems can be be found in \cite{Fazio:2019:NIT}.

Moreover, T\"opfer's method has been extended to classes of problems in boundary layer theory involving one or more physical parameters.
This kind of extension was considered first by Na \cite{Na:1970:IVM}, see also NA \cite[Chapters 8-9]{Na:1979:CME} for an extensive survey on this subject.

Finally, an iterative extension of the transformation method has been introduced, for the numerical solution of free BVPs, by Fazio \cite{Fazio:1990:SNA}. 
This iterative extension has been applied to several problems of interest: free boundary problems \cite{Fazio:1990:SNA,Fazio:1997:NTE,Fazio:1998:SAN},
a moving boundary hyperbolic problem \cite{Fazio:1992:MBH}, Homann and Hiemenz problems governed by the Falkner-Skan equation in \cite{Fazio:1994:FSE},
one-dimensional parabolic moving boundary problems \cite{Fazio:2001:ITM}, two variants of the Blasius problem \cite{Fazio:2009:NTM}, namely: a boundary layer problem over moving surfaces, studied first by Klemp and Acrivos \cite{Klemp:1972:MBL}, and a boundary layer problem with slip boundary condition, that has found application in the study of gas and liquid flows at the micro-scale regime \cite{Gad-el-Hak:1999:FMM,Martin:2001:BBL}, parabolic problems on unbounded domains \cite{Fazio:2010:MBF} and, recently, see \cite{Fazio:2015:ITM}, a further variant of the Blasius problem in boundary layer theory: the so-called Sakiadis problem \cite{Sakiadis:1961:BLBa,Sakiadis:1961:BLBb}.
A recent review dealing with, the derivation and application of, ITM can be be found, by the interested reader, in \cite{Fazio:2019:ITM}.
A unifying framework, providing proof that the non-ITM is a special instance of the ITM and consequently can be derived from it, has been the argument of the paper \cite{Fazio:2020:SIT}.

\section{Blasius problem with extended boundary conditions}
The Blasius problem with extended boundary conditions is given by, see White \cite{White:2006:VFF}, Klemp and Acrivos \cite{Klemp:1972:MIB} and Fang and Lee \cite{Fang:2005:MWB},
\begin{align}\label{eq:ExBlasius} 
& {\displaystyle \frac{d^{3}f}{d\eta^3}} + f
{\displaystyle \frac{d^{2}f}{d\eta^2}} = 0 \nonumber \\[-1ex]
&\\[-1ex]
& f(0) = 0 \ , \qquad {\displaystyle \frac{df}{d\eta}}(0) = P_1 + P_2 \; {\displaystyle \frac{d^2f}{d\eta^2}}(0) \ , \qquad {\displaystyle \frac{df}{d\eta}}(\eta) \rightarrow 1 \quad \mbox{as}
\quad \eta \rightarrow \infty \ , \nonumber 
\end{align}
where $P_1 = \frac{U_w}{U_{\infty}}$, which for $U_w > 0$ can be positive with the same direction as the free stream velocity and for $U_w < 0$ is negative opposite to the free stram velocity, and $P_2 = \frac{U_{slip}}{U_{\infty}} = \left(\frac{2}{P_3}-1\right) K_{}n,x Re_x^{1/2}$ is a dimensionless parameter with $K_{n,x} = \frac{1}{x}$, and $Re_x = \frac{U_{\infty}x}{2\nu}$. 
We notice here, that the problem (\ref{eq:ExBlasius}) when $P_1 = P_2 = 0$ reduces to the celebrated Blasius problem.

\subsection{The non-ITM} 
In this section, we assume that we need to find the behaviour of the missing initial condition with respect to the variation of the values of the involved parameters, that is $P_1$ and $P_2$ should get several different values but these values are no fixed in advance.
The applicability of a non-ITM to the Blasius problem is a consequence of both: the invariance of the governing differential equation and the two boundary conditions at $\eta = 0$, and the non-invariance of the asymptotic boundary condition, as $\eta$ goes to infinity, under the scaling group of point transformation.
In order to apply a non-ITM to the BVP (\ref{eq:ExBlasius}) we investigate its invariance with respect to the extended scaling group
\begin{equation}\label{eq:scaling}
f^* = \lambda f \ , \qquad \eta^* = \lambda^{-1} \eta \ , \qquad P_1^* = \lambda^{\delta_1} P_1 \ , \qquad P_2^* = \lambda^{\delta_2} P_2 \ .   
\end{equation}
We find that the Blasius problem with extended boundary conditions (\ref{eq:ExBlasius}) is invariant under (\ref{eq:scaling}) iff
\begin{equation}\label{eq:scaling:condition}
\delta_1 = 2 \ , \qquad \delta_2 = -1 \ .
\end{equation}
Now, we can integrate the Blasius equation in (\ref{eq:ExBlasius}) written in the starred variables on $[0, \eta^*_\infty]$, where $\eta^*_\infty$ is a suitable truncated boundary, with initial conditions
\begin{equation}\label{eq:ICs}
f^*(0) = 0 \ , \qquad \frac{df^*}{d\eta^*}(0) = P_1^*+P_2^* \; {\displaystyle \frac{d^2f^*}{d\eta^{*2}}(0)} \ , \quad {\displaystyle \frac{d^2f^*}{d\eta^{*2}}(0)} = 1 \ ,
\end{equation}
in order to compute an approximation $\frac{df^*}{d\eta^*}(\eta^*_\infty)$ for $\frac{df^*}{d\eta^*}(\infty)$ and the corresponding value of $\lambda$ according to the equation
\begin{equation}\label{eq:lambda}
\lambda = \left[ \frac{d f^*}{d \eta^{*}}(\eta^*_\infty) \right]^{1/2} \ .   
\end{equation} 
Once the value of $\lambda$ has been computed by equation (\ref{eq:lambda}), we can find the missed initial condition by the equation
\begin{equation}\label{eq:MIC}
\frac{d^2f}{d\eta^{2}}(0) =  \lambda^{2\delta-1}\frac{d^2f^*}{d\eta^{*2}}(0) \ ,
\end{equation}
and the values of $P_1$ and $P_2$ by the relations
\begin{equation}\label{eq:Parameters}
P_1 = \lambda_1^{-2} \ , \qquad P_2 = \lambda_2 \ .
\end{equation}
Moreover, the numerical solution of the original BVP (\ref{eq:ExBlasius}) can be computed by rescaling the numerical solution of the IVP.
In this way, we get the solution of a given BVP by solving a related IVP.

\subsection{The ITM} 
In this section, we assume that we need to compute the numerical solution for given values of the involved parameters, that is $P_1$ and $P_2$ are now fixed.
We need now to consider the invariance of the initial conditions with respect to the extended scaling group of point transformations
\begin{equation}\label{eq:scaling2}
f^* = \lambda f \ , \qquad \eta^* = \lambda^{-1} \eta \ , \qquad h^* = \lambda^{\sigma} h \ .   
\end{equation}
This new scaling group involves the scaling of the fictitious parameter $h$ that will be used to force the initial conditions to be invariant. 
Now, we can integrate the Blasius equation in (\ref{eq:ExBlasius}) written in the star variables on $[0, \eta^*_\infty]$, where $\eta^*_\infty$ is a suitable truncated boundary, with initial conditions
\begin{equation}\label{eq:ICs2}
f^*(0) = 0 \ , \qquad \frac{df^*}{d\eta^*}(0) = h^{2/\sigma} P_1 + h^{-1/\sigma} P_2 \; {\displaystyle \frac{d^2f^*}{d\eta^{*2}}(0)} \ , \quad {\displaystyle \frac{d^2f^*}{d\eta^{*2}}(0)} = 1 \ ,
\end{equation}
in order to compute an approximation $\frac{df^*}{d\eta^*}(\eta^*_{\infty})$ for $\frac{df^*}{d\eta^*}(\infty)$ and the corresponding value of $\lambda$ again by equation (\ref{eq:lambda})
Once the value of $\lambda$ has been computed by equation (\ref{eq:lambda}), we can find the missed initial condition again from equation (\ref{eq:MIC}).
In the ITM we proceed as follows: we set the values of $P_1$, ${P_2}$, $h^*$, $\sigma$ and $\eta_{\infty}^*$ and integrate the IVP on $[0, \eta_{\infty}^*]$.
Naturally, choosing $h^*$ arbitrarily we do not obtain the value $h = 1$, however, we can apply a root-finder method, like bisection, secant, regula-falsi, Newton or quasi-Newton root-finder, because the required value of $h$ can be considered as the root of the implicit defined, transformation function
\begin{equation}\label{eq:Tfunction}
\Gamma(h^*) = \lambda^{-\sigma} h -1  \ .   
\end{equation}

Of course, any positive value of $\sigma$ can be chosen, and in the following, for the sake of simplicity, we set $\sigma = 10$.
Moreover, as a termination criterion for our root-finder we used $|\Gamma(h^*)| < Tol$ with $Tol = 10^{-5}$.

\section{Numerical results}
In this section, we report the numerical results computed with our non-ITM and ITM. 
To compute the numerical solution, we used the eighth order Runge-Kutta method \cite[p. 180]{Butcher:NMO:2003} with constant step size.

First of all, we start with the results obtained by the non-ITM
In table \ref{tab:NITM:missingIC} we report the chosen parameter values, the computed values of the involved parameters and the missing initial condition $\frac{d^2f}{d\eta^2}(0)$.
\begin{table}[!hbt]
\caption{Numerical data and results.}
\vspace{.5cm}
\renewcommand\arraystretch{1.3}
	\centering
		\begin{tabular}{llllr@{.}l}
\hline 
{$P_1^*$} & {$P_2^*$} & {$P_1$} & {$P_2$} & \multicolumn{2}{c}%
{$ {\displaystyle \frac{d^2f}{d\eta^2}(0)}$} \\[1.2ex]
\hline
0.25 & 0.25 & 0.140225769 &  0.333807506 & 0 & 42007973468 \\  
0.5  & 0.5  & 0.241979004 &  0.336675506 & 0 & 33667550559 \\  
0.75 & 0.75 & 0.309184205 &  1.168108665 & 0 & 26468787856 \\  
1    & 1    & 0.353764405 &  1.681291175 & 0 & 21041233684 \\  
1.5  & 1.5  & 0.405947260 &  2.883381325 & 0 & 14078861396 \\  
2    & 2    & 0.433836425 &  4.294197226 & 0 & 10102852811 \\  
2.5  & 2.5  & 0.450478633 &  5.889425257 & 0 & 07648940496 \\  
5    & 5    & 0.481068451 & 16.119500068 & 0 & 02984388156 \\   
\hline			
		\end{tabular}
	\label{tab:NITM:missingIC}
\end{table}
As it is easily seen from the results listed in table \ref{tab:NITM:missingIC} we are not in the position to plot the data by fixing one of the two parameters, usually $P_1$, and plotting the missing initial condition versus the other parameter. 
Of course, this is a drawback of our non-ITM.
However, when we are required to do just these kinds of plots we can apply the described ITM.

We report now, the numerical results obtained by the ITM.
As a root-finder we applied the simple bisection method with the termination criterion $|\Gamma(h^*)| < Tol$ with $Tol = 10^{-5}$.
In table \ref{tab:ITM:ITERA} we report a sample iteration of the bisection method.
\begin{table}[!hbt]
\caption{Bisection method iterations for $P_1 = 0.5 $ and $P_2 = 0$.}
\vspace{.5cm}
\renewcommand\arraystretch{1.3}
	\centering
		\begin{tabular}{llr@{.}l}
\hline 
{$h^*$} & {$\lambda$} & \multicolumn{2}{c}%
{$\Gamma(h^*)$} \\[1.2ex]
\hline
0.75             &             & $-$0 & 424804078 \\ 
1.75             &             &    0 & 118076477 \\
1.25             & 1.389163618 & $-$0 & 100177989 \\ 
1.5              & 1.466575876 &    0 & 022790586 \\
1.375            & 1.425023536 & $-$0 & 035103656 \\
1.4375           & 1.445108710 & $-$0 & 005265147 \\
1.46875          & 1.455672550 &    0 & 008983786 \\
1.453125         & 1.450347802 &    0 & 001914850 \\
1.4453125        & 1.447717501 & $-$0 & 001661237 \\
1.44921875       & 1.449029969  &   0 & 000130281 \\
1.447265625      & 1.448373064 & $-$0 & 0007646088 \\
1.4482421875     & 1.448701349 & $-$0 & 0003169467 \\
1.44873046875    & 1.448865617 & $-$0 & 0000932785 \\
1.448974609375   & 1.448947782 &    0 & 0000185148 \\
1.4488525390625  & 1.448906697 & $-$0 & 0000373785 \\
1.44891357421875 & 1.448927239 & $-$0 & 0000094310 \\
\hline			
		\end{tabular}
	\label{tab:ITM:ITERA}
\end{table}

\begin{figure}[!hbt]
	\centering
\psfragscanon 
\psfrag{f2}[c][c]{$\frac{d^2f}{d\eta^2}$}  
\psfrag{P1}[c][c]{$P_1$}  
\psfrag{P20}[c][c]{$P2 = 0$} 
\includegraphics[width=13cm,height=13cm]{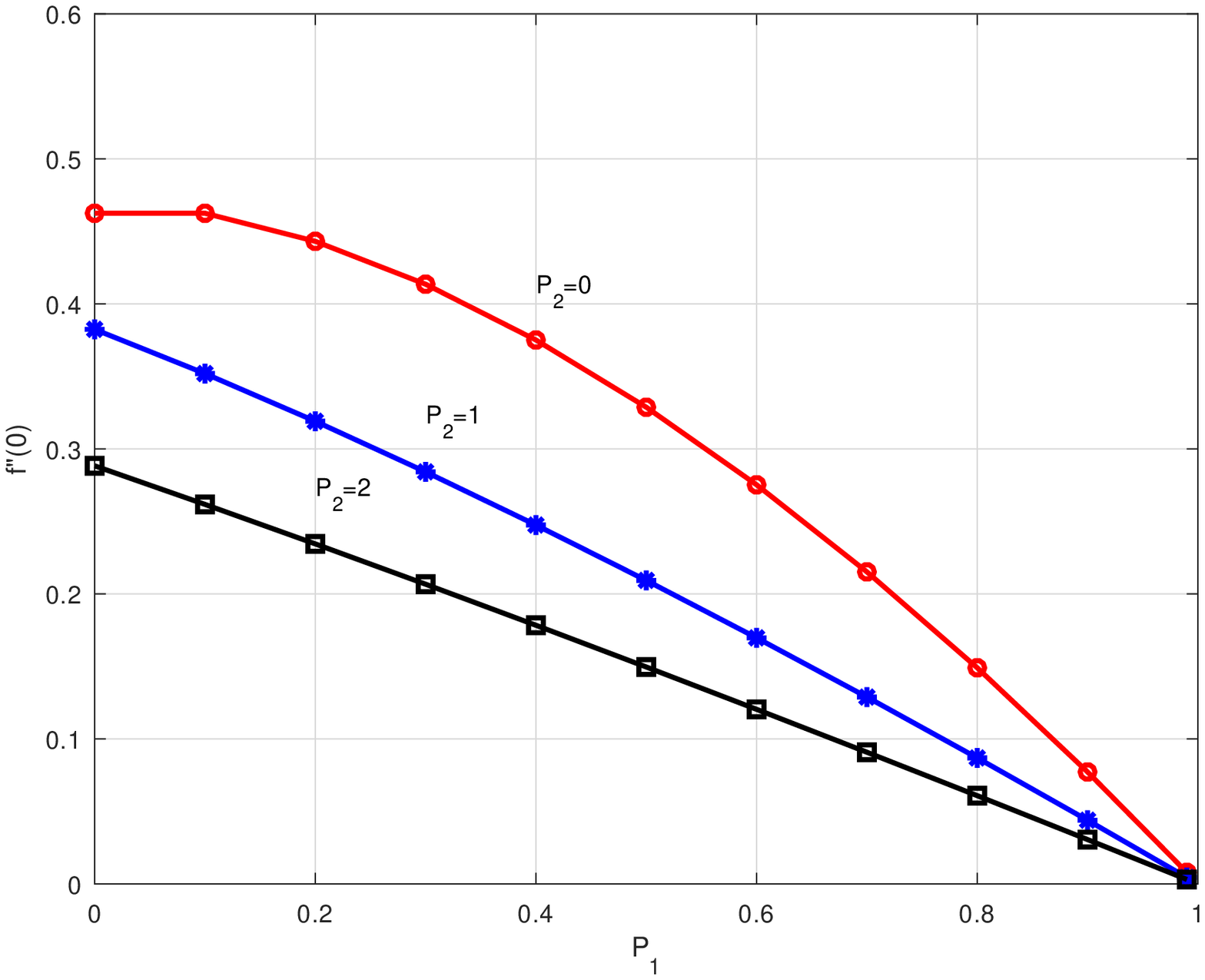}
\caption{Numerical results of the missing initial condition versus $P_1$, here $P_2 = 0, 1, 2$.}
	\label{fig:f2P1}
\end{figure}

In figure \ref{fig:f2P1} shows the behaviour of the missing initial condition versus $P_1$ with three values of the other parameter, namely $P_2 = 0, 1, 2$.

\begin{figure}[!hbt]
	\centering
\psfrag{e}[1][]{$\eta, \eta^*$}  
\psfrag{fd}[l][]{$\ds \frac{df}{d\eta}$} 
\psfrag{ddf}[l][]{$\ds \frac{d^2f}{d\eta^2}$} 
\psfrag{fd*}[l][]{$\ds \frac{df^*}{d\eta^*}$} 
\psfrag{ddf*}[l][]{$\ds \frac{d^2f^*}{d{\eta^*}^2}$} 
\includegraphics[width=13cm,height=13cm]{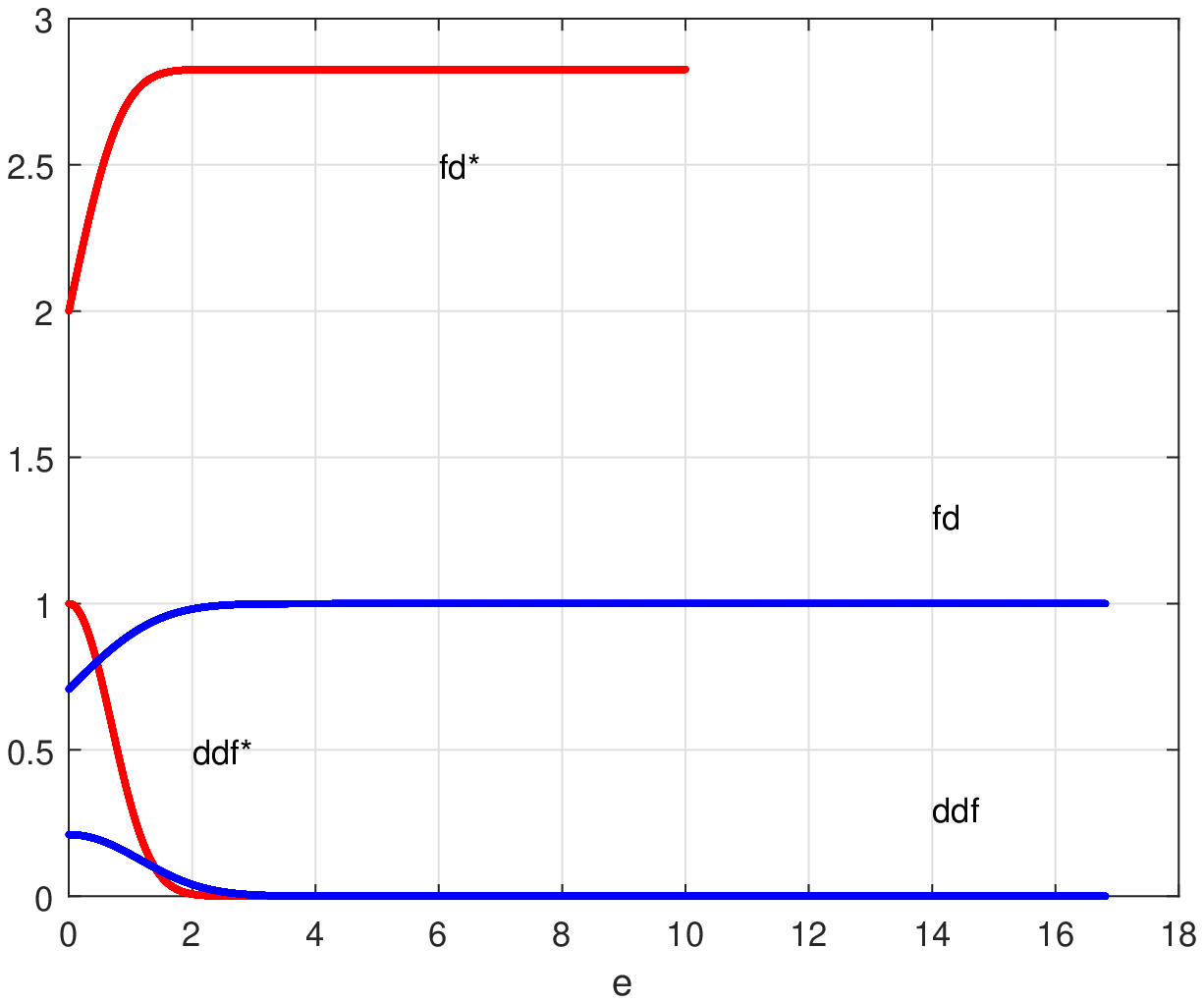}
\caption{Numerical results of the non-ITM for (\ref{eq:ExBlasius}) with $P_1 = P_2 = 1$. The starred variables problem and the original problem solution components found after rescaling.}
	\label{fig:BlasiusExBCs}
\end{figure}

As an example, figure \ref{fig:BlasiusExBCs} shows the solution of the Blasius problem with extended boundary condition in the particular case when we set $P_1^* = P_2^* = $1.
For the results shown in this figure we used $\Delta \eta = 0.001$ and $\eta^*_{\infty} = 10$.
Let us notice here that, by rescaling, we get $\eta^*_{\infty} < \eta_{\infty}$, and this is convenient for the user because it means that we need to make less computational effort to get the wanted numerical solution. 

As mentioned before, the case $P_1 = P_2 = 0$ is the Blasius problem.
In this particular case, our non-ITM reduces to the original method defined by T\"opfer \cite{Topfer:1912:BAB}, and the computed skin friction coefficient value, namely $0.469599988361$, obtained with $\Delta \eta = 0.0001$ and $\eta^*_{\infty} = 10$, is in very good agreement with the values available in the literature, see for instance the value $0.469599988361$ computed by Fazio \cite{Fazio:1992:BPF} by a free boundary formulation of the Blasius problem.

\section{Concluding remarks}
The main contribution of this paper is the extension of the non-ITM, proposed by T\"opfer \cite{Topfer:1912:BAB} and defined for the numerical solution of the celebrated Blasius problem \cite{Blasius:1908:GFK}, to a Blasius problem with extended boundary conditions.
This method, that makes use of the invariance of two physical parameters, allows us to solve numerically the Blasius problem with extended boundary conditions by solving a related IVP and then rescaling the obtained numerical solution. 
However, in this way, we cannot fix in advance the physical parameters, and if we need just to compute the numerical solution for given values of the two parameters we have to define an iterative extension of our TM.

\vspace{1.5cm}

\noindent {\bf Acknowledgement.} {The research of this work was 
partially supported by the University of Messina and by the GNCS of INDAM.}


\end{document}